\documentclass[12pt]{article}

\usepackage{amsmath}
\usepackage{amssymb}
\usepackage{amscd}
\usepackage{amsthm}
\usepackage{indentfirst}
\usepackage[english,frenchb]{babel}

\newtheorem{thm}{Th\'eor\`eme}[section]
\newtheorem{lem}[thm]{Lemme}
\newtheorem{cor}[thm]{Corollaire}
\newtheorem{prop}[thm]{Proposition}
\newtheorem{dfn}[thm]{D\'efinition}

\newtheorem{question}[thm]{Question}
\theoremstyle{remark}
\newtheorem{rmq}[thm]{Remarque}
\newtheorem{exm}[thm]{Exemple}
\newtheorem{ntt}[thm]{Notations}

\DeclareMathOperator{\Spec}{Spec}

\DeclareMathOperator{\pic}{Pic}
\DeclareMathOperator{\ord}{ord}

\DeclareMathOperator{\rpic}{RPic}
\DeclareMathOperator{\rk}{rang}

\newcommand{\A}{\mathcal{A}}
\newcommand{\B}{\mathcal{B}}
\newcommand{\C}{\mathcal{C}}

\newcommand{\X}{\mathcal{X}}

\newcommand{\E}{\mathcal{E}}

\newcommand{\gm}{\mathbf{G}_{{\rm m}}}
\newcommand{\gmk}{\mathbf{G}_{{\rm m},K}}
\newcommand{\ga}{\mathbf{G}_{{\rm a}}}

\newcommand{\psibemol}{\psi^{\flat}}

\newcommand{\znz}{(\mathbb{Z}/n\mathbb{Z})}
\newcommand{\zpnz}{(\mathbb{Z}/p^n\mathbb{Z})}

\newcommand{\zpz}{(\mathbb{Z}/p\mathbb{Z})}

\newcommand{\Hom}{\underline{\rm Hom}}

\newcommand{\Ext}{\underline{\rm Ext}^1}
\newcommand{\ext}{{\rm Ext}^1}

\begin{document}

\title{Invariants de classes : exemples de non-annulation en dimension sup\'erieure}

\author{Jean Gillibert}

\date{8 d\'ecembre 2006}

\maketitle

\begin{abstract}
Le class-invariant homomorphism permet de mesurer la structure galoisienne des torseurs --- sous un sch\'ema en groupes fini et plat $G$ --- qui sont dans l'image du cobord associ\'e \`a une isog\'enie, de noyau $G$, entre des (mod\`eles de N\'eron de) vari\'et\'es ab\'eliennes. Quand les vari\'et\'es sont des courbes elliptiques \`a r\'eduction semi-stable et que l'ordre de $G$ est premier \`a $6$, on sait que cet homomorphisme s'annule sur les points de torsion. Dans cet article, en nous servant de restrictions de Weil de courbes elliptiques, nous construisons, pour tout nombre premier $p>2$, une vari\'et\'e ab\'elienne $A$ de dimension $p$ munie d'une isog\'enie (de noyau $\mu_p$) dont le cobord est surjectif. Si $A$ est de rang nul, et si la $p$-partie du groupe de Picard de la base est non triviale, nous obtenons ainsi un exemple o\`u le class-invariant homomorphism ne s'annule pas sur les points de torsion.

\medskip

\begin{otherlanguage}{english}
\begin{center}
\textbf{Abstract}
\end{center}

The so-called class-invariant homomorphism $\psi$ measures the Galois module structure of torsors---under a finite flat group scheme $G$---which lie in the image of a coboundary map associated to an isogeny between (N\'eron models of) abelian varieties with kernel $G$. When the varieties are elliptic curves with semi-stable reduction and the order of $G$ is coprime to $6$, it is known that the homomorphism $\psi$ vanishes on torsion points. In this paper, using Weil restrictions of elliptic curves, we give the construction, for any prime number $p>2$, of an abelian variety $A$ of dimension $p$ endowed with an isogeny (with kernel $\mu_p$) whose coboundary map is surjective. In the case when $A$ has rank zero and the $p$-part of the Picard group of the base is non-trivial, we obtain examples where $\psi$ does not vanish on torsion points.
\end{otherlanguage}
\end{abstract}




\section{Introduction}

Soit $S$ un sch\'ema, et soit $G$ un $S$-sch\'ema en groupes commutatif, fini localement libre. Soit $H^1(S,G)$ le premier groupe de cohomologie de $S$ \`a valeurs dans $G$, calcul\'e pour la topologie fppf sur $S$. On sait que les \'el\'ements de $H^1(S,G)$ s'identifient aux (classes d'isomorphie de) $G$-torseurs sur $S$, en particulier quand $G=\Gamma_S$ est un $S$-groupe constant fini les \'el\'ements de $H^1(S,\Gamma_S)$ sont les rev\^etements galoisiens non ramifi\'es de groupe $\Gamma$.

Une fa\c con naturelle de construire des $G$-torseurs est de consid\'erer une suite exacte de $S$-sch\'emas en groupes (pour la topologie fppf sur $S$) de noyau $G$,
\begin{equation}
\label{unesuite}
\begin{CD}
0 @>>> G @>>> P @>f>> Q @>>> 0 \\
\end{CD}
\end{equation}
qui donne lieu \`a un morphisme cobord
\begin{equation}
\label{cobord}
\begin{CD}
\Delta:Q(S) @>>> H^1(S,G). \\
\end{CD}
\end{equation}
De fa\c con imag\'ee, \`a toute section $x\in Q(S)$, on associe un $G$-torseur qui repr\'esente le faisceau \og image r\'eciproque de $x$ par $f$ \fg{}. Les exemples de telles situations abondent en g\'eom\'etrie alg\'ebrique. Citons par exemple la suite exacte de Kummer, dans laquelle $G=\mu_n$ et $f:\gm\rightarrow\gm$ est le morphisme d'\'el\'evation \`a la puissance $n$-i\`eme dans le groupe multiplicatif $\gm$. Dans le cas o\`u $S$ est de caract\'eristique $p$, citons \'egalement la suite exacte d'Artin-Schreier, dans laquelle $G=\zpz_S$ et $f:\ga\rightarrow\ga$ est le morphisme $\wp={\rm id}-F$ dans le groupe additif $\ga$, o\`u $F$ d\'esigne le morphisme de Frobenius.

Une question naturelle se pose : \'etant donn\'e $G$, peut-on trouver une suite exacte de la forme (\ref{unesuite}) telle que le cobord $\Delta$ soit surjectif, et que peut-on imposer sur $P$ et $Q$ ?

Consid\'erons le cas particulier o\`u $S=\Spec(R)$, avec $R$ un anneau de Dedekind de type fini sur $\mathbb{Z}$, dont le corps de fractions $K$ est un corps de nombres. Si $G=\mu_n$ alors dans \cite{gil3} on construit une suite exacte de la forme (\ref{unesuite}) dans laquelle $P$ (donc \'egalement $Q$) est un $S$-tore, et dont le cobord est surjectif. Par contre, si $G=\znz_S$ avec $n$ non inversible sur $S$, alors il est impossible de plonger $G$ dans un $S$-tore.

D'apr\`es un r\'esultat de Raynaud (voir \cite[Th\'eor\`eme 3.1.1]{bbm}) il est possible, localement pour la topologie de Zariski sur $S$, de plonger $G$ dans un $S$-sch\'ema ab\'elien (\emph{i.e.} un $S$-sch\'ema en groupes lisse, propre, \`a fibres connexes). Cependant, le probl\`eme n'a pas toujours de solution au niveau global : si $R=\mathbb{Z}$ alors, d'apr\`es un c\'el\`ebre r\'esultat de Fontaine \cite{font}, il n'existe pas de sch\'ema ab\'elien sur $S$. Il est donc indispensable d'affaiblir les exigences en autorisant de la mauvaise r\'eduction. Nous arrivons donc \`a la question suivante :
existe-t-il une suite de la forme (\ref{unesuite}) dans laquelle $P$ est le mod\`ele de N\'eron d'une $K$-vari\'et\'e ab\'elienne, et dont le cobord est surjectif ?
En adaptant au cas pr\'esent la strat\'egie utilis\'ee dans la preuve de la proposition 3.4 de \cite{gil3},
nous montrons ici le r\'esultat suivant (cf. le paragraphe \ref{surjectif}).

\begin{thm}
\label{sgsemistable}
Supposons que $G$ soit un sous-groupe fini et plat du mod\`ele de N\'eron $\X$ d'une $K$-vari\'et\'e ab\'elienne ayant r\'eduction semi-stable en les points de $S$ de caract\'eristique divisant l'ordre de $G$. Alors il existe une $K$-vari\'et\'e ab\'elienne dont le mod\`ele de N\'eron $\A$ contient $G$ comme sous-groupe, et telle que le cobord associ\'e \`a la suite exacte
$$
\begin{CD}
0 @>>> G @>>> \A @>>> B @>>> 0 \\
\end{CD}
$$
soit surjectif.
\end{thm}

La vari\'et\'e ab\'elienne $A_K$ est en fait de la forme $\Re_{K'/K}(X_{K'})$, o\`u $K'$ est une extension finie de $K$ telle que la fl\`eche de changement de base $H^1(S,G)\rightarrow H^1(S',G_{S'})$ soit nulle, $S'$ \'etant le normalis\'e de $S$ dans $K'$, et o\`u $\Re_{K'/K}$ d\'esigne la restriction des scalaires \`a la Weil (voir \cite[\S\/7.6]{ray}). En outre, le quotient $B:=\A/G$ est repr\'esentable par un $S$-sch\'ema en groupes lisse, s\'epar\'e, et de type fini, dont la fibre g\'en\'erique est une vari\'et\'e ab\'elienne $B_K$ (voir les lemmes \ref{lem22} et \ref{lem24}).

La motivation initiale qui nous a men\'e \`a ces questions est l'\'etude de la conjecture de M.~J.~Taylor sur le \emph{class-invariant homomorphism} pour des vari\'et\'es ab\'eliennes de dimension strictement sup\'erieure \`a $1$. Dans cet article, nous d\'efinissons cet homomorphisme (que nous noterons $\psi$ et que nous appellerons homomorphisme de classes) comme \'etant la compos\'ee du cobord $\Delta$ de (\ref{cobord}) avec un morphisme
\[
\begin{CD}
\pi:H^1(S,G) @>>> \pic(G^D)\\
\end{CD}
\]
d\'efini par Waterhouse \cite[Theorem $5$]{w}, et qui s'interpr\`ete en termes de structure galoisienne de torseurs. Plus pr\'ecis\'ement, si $Z$ est un $G$-torseur sur $S$, alors $Z$ est affine, et l'action de $G$ sur $Z$ donne naissance \`a une action de l'alg\`ebre du groupe $G^D$ (\emph{i.e.} l'alg\`ebre duale de l'alg\`ebre du groupe $G$) sur l'alg\`ebre de $Z$. L'homomorphisme $\pi$ mesure cette structure de module. Nous noterons $\rpic(G^D)$ l'image de $\pi$, et nous l'appellerons le groupe des classes r\'ealisables.

D'apr\`es ce qui pr\'ec\`ede, on peut traduire les r\'esultats ci-dessus en termes de surjectivit\'e de l'homomorphisme de classes sur le groupe des classes r\'ealisables.

Dans \cite{t2}, M.~J.~Taylor a conjectur\'e l'annulation de $\psi$ sur les points de torsion, dans le cas o\`u $P$ est une $S$-courbe elliptique. Si l'ordre de $G$ est premier \`a $6$, ce r\'esultat est connu. D'abord d\'emontr\'e par Srivastav et Taylor \cite{st} pour une courbe \`a multiplication complexe, puis par Agboola \cite{a2} sans l'hypoth\`ese de multiplication complexe, et enfin par Pappas \cite{p1} en rempla\c cant $S$ par un sch\'ema quelconque, il a \'egalement \'et\'e g\'en\'eralis\'e par l'auteur (voir \cite{gil1} et \cite{gil2}) dans le cas o\`u $P$ est un $S$-sch\'ema en groupes semi-stable dont la fibre g\'en\'erique est une $K$-courbe elliptique.

Si $P$ est une $S$-courbe elliptique, et si l'ordre de $G$ est une puissance de $2$, alors il existe des exemples de non-annulation de $\psi$ sur des points de torsion, d\^us \`a W. Bley et M. Klebel \cite{bk}, ainsi qu'\`a Ph. Cassou-Nogu\`es et A. Jehanne \cite{cnj} dans le cas d'une courbe elliptique \`a multiplication complexe.

Dans \cite[Theorem C]{p1}, pour tout choix de deux nombres premiers $r\neq\ell$, on construit une courbe affine lisse $D$ sur un corps fini de caract\'eristique $r$, et un $D$-sch\'ema ab\'elien de dimension $2$ muni d'un point de $\ell$-torsion sur lequel l'homomorphisme de classes associ\'e \`a la multiplication par $\ell$ ne s'annule pas. Cependant, aucun r\'esultat semblable n'est connu dans le cas o\`u le sch\'ema de base est le spectre de l'anneau des entiers d'un corps de nombres.

Notre r\'esultat donne une nouvelle m\'ethode pour construire des points de torsion qui ne sont pas dans le noyau de $\psi$. Soit $p$ un nombre premier, consid\'erons le cas o\`u $G=\mu_p$, alors $G^D=\zpz_S$ et $\rpic(G^D)=\pic(S)[p]$, donc on obtient un homomorphisme de classes $\psi$ surjectif sur $\pic(S)[p]$. Pour obtenir un contre-exemple \`a la conjecture de Taylor, il suffit de se placer dans le cas o\`u $\pic(S)[p]$ n'est pas nul, et $B(S)$ est de torsion. Plus pr\'ecis\'ement, nous obtenons le r\'esultat suivant (cf. le th\'eor\`eme \ref{bigthm}).

\begin{thm}
\label{undeux}
Soient
\begin{enumerate}
\item $p\geq 3$ un nombre premier,
\item $K$ un corps quadratique imaginaire tel que $\pic(\mathcal{O}_K)[p]\neq 0$ ; si $p=3$ on exige de plus que $p$ soit non ramifi\'e dans $K$,
\item $K'$ une extension de $K$ dans laquelle les \'el\'ements de $\pic(\mathcal{O}_K)[p]$ capitulent,
\item $E_{K'}$ une $K'$-courbe elliptique, ayant bonne r\'eduction au-dessus de $p$, telle que $E_{K'}(K')$ soit un groupe fini contenant un \'el\'ement d'ordre $p$.
\end{enumerate}
Soit $S=\Spec(\mathcal{O}_K)$. Alors il existe une $K$-vari\'et\'e ab\'elienne de dimension $[K':K]$, de mod\`ele de N\'eron $\A$, et une suite exacte
$$
\begin{CD}
0 @>>> \mu_p @>>> \A @>>> B @>>> 0 \\
\end{CD}
$$
o\`u $B(S)$ est un groupe de torsion, telle que l'homomorphisme de classes associ\'e
$$
\begin{CD}
\psi:B(S) @>>> H^1(S,\mu_p) @>>> \rpic(\zpz_S)\simeq\pic(\mathcal{O}_K)[p] \\
\end{CD}
$$
soit surjectif, donc non nul sur les points de torsion.
\end{thm}

Dans la pratique, la finitude de $E_{K'}(K')$ est assez difficile \`a v\'erifier par le calcul. Nous avons pour l'instant des exemples pour $p=3$ et $p=5$. Il est cependant raisonnable de croire que de tels exemples existent pour d'autres valeurs de $p$.

Enfin, sous certaines hypoth\`eses, nous donnons une interpr\'etation de $\psi$ en termes de la restriction de certains $\gm$-torseurs prolongeant le fibr\'e de Poincar\'e sur $B_K\times_K B_K^t$, o\`u $B_K^t$ d\'esigne la vari\'et\'e duale de $B_K$. Plus pr\'ecis\'ement, soit $\B$ le mod\`ele de N\'eron de $B_K$, et soit $\B^{\Lambda}$ (voir les notations \ref{notation1}) le plus petit sous-groupe ouvert de $\B$ \`a travers lequel le morphisme canonique $B\rightarrow\B$ se factorise. Supposons qu'il existe un morphisme $i:G^D\rightarrow\B^{t,\Lambda'}$ prolongeant l'inclusion canonique $G_K^D\subseteq B_K^t$, o\`u $\Lambda'$ d\'esigne l'orthogonal de $\Lambda$ sous l'accouplement de monodromie (cf. le paragraphe 2.2), ce qui est v\'erifi\'e par exemple si $G^D$ est \'etale sur $S$. On peut alors consid\'erer le morphisme $\psibemol$ donn\'e par
\begin{equation}
\begin{split}
\psibemol:\B^{\Lambda}(S)\;\longrightarrow & \;\rpic(G^D)\\
(x:S\rightarrow \B^{\Lambda})\;\longmapsto & \;(x\times i)^*(t(W'))\\
\end{split}
\end{equation}
o\`u $W'$ est l'unique biextension de $(\B^{\Lambda}, \B^{t,\Lambda'})$ par $\gm$ qui prolonge la biextension de Weil, et o\`u $t(W')$ est le $\gm$-torseur associ\'e \`a $W'$. On montre alors que $\psi$ est la compos\'ee de $\psibemol$ avec la fl\`eche $B(S)\rightarrow\B^{\Lambda}(S)$, ce qui g\'en\'eralise la description de $\psi$ donn\'ee par Agboola \cite[Theorem 1]{a1} dans le cas des sch\'emas ab\'eliens.
Nous sommes alors en mesure d'en d\'eduire l'\'enonc\'e qui suit (cf. le corollaire \ref{petitcor}).

\begin{thm}
\label{untrois}
Soit $N>0$ un entier. Il existe une $K$-vari\'et\'e ab\'elienne $B_K$, un sous-groupe $\Lambda$ du groupe des composantes du mod\`ele de N\'eron $\B$ de $B_K$, et un point de $N$-torsion $y\in \B^{t,\Lambda'}(S)$, tels que l'application
\begin{equation*}
\begin{split}
\B^{\Lambda}(S)\;\longrightarrow & \;\pic(S)[N]\\
(x:S\rightarrow \B^{\Lambda})\;\longmapsto & \;(x\times y)^*(t(W'))\\
\end{split}
\end{equation*}
soit surjective,
o\`u $W'$ est l'unique biextension de $(\B^{\Lambda}, \B^{t,\Lambda'})$ par $\gm$ qui prolonge la biextension de Weil, et o\`u $t(W')$ est le $\gm$-torseur associ\'e \`a $W'$.
\end{thm}

Pour conclure, nous proposons de formuler la question suivante, qui contient la conjecture de Taylor en dimension $1$, et qui inclut le cas des corps de fonctions. Notons ici que les contre-exemples de Pappas n'excluent pas une r\'eponse positive \`a cette question.

\begin{question}
Soit $S$ le spectre d'un anneau de Dedekind dont le corps de fraction $K$ est un corps global (\emph{i.e.} un corps de nombres, ou une extension de degr\'e de transcendance un d'un corps fini), et soit $d$ un entier. Existe-t-il un entier $Q(K,d)$ tel que, pour toute $K$-vari\'et\'e ab\'elienne $C_K$ de dimension $d$, de mod\`ele de N\'eron $\C$ \`a r\'eduction semi-stable sur $S$, la restriction de $t(W)$ \`a $\C^{\circ}[m]\times\C^t[m]$ soit triviale pour $m$ premier \`a $Q(K,d)$ ?
\end{question}

\begin{rmq}
(1) D'apr\`es \cite[th\'eor\`eme 4.1]{gil1}, $Q(K,1)=6$ convient pour tout $K$. De plus, d'apr\`es \cite[Theorem C]{p1}, pour tout $p$ premier, il existe un corps de fonctions $K$, de caract\'eristique diff\'erente de $p$, tel que $p|Q(K,d)$ pour tout entier $d\geq 2$. Ainsi l'entier $d=1$ est le seul pour lequel il existe une constante $Q(K,d)$ qui est ind\'ependante de $K$.

(2) Si $K$ est un corps de nombres et si $p$ un nombre premier, satisfaisant les hypoth\`eses du th\'eor\`eme \ref{undeux}, alors $p|Q(K,p)$.

(3) Cette question est \`a mettre en parall\`ele avec la conjecture 1.8 de Pappas dans \cite{p3}. Il existe en effet un lien entre la caract\'eristique d'Euler \'equivariante consid\'er\'ee par Pappas, et les valeurs de l'homomorphisme de classes. Pour plus de d\'etails, nous renvoyons le lecteur \`a \cite{p2}.
\end{rmq}


\section{Sous-groupes finis}

Dans tout le texte, $K$ est un corps de nombres, $R$ est un anneau de Dedekind (de type fini en tant que $\mathbb{Z}$-alg\`ebre) de corps de fractions $K$, et $S=\Spec(R)$.
Tous les sch\'emas en groupes que nous consid\`ererons ici sont commutatifs par convention.

Dans cette section, nous fixons une $K$-vari\'et\'e ab\'elienne $A_K$, et nous notons $\A$ son mod\`ele de N\'eron sur $S$.  Si $G$ est un $S$-sch\'ema en groupes (commutatif) fini et plat, nous noterons $\ord(G)$ l'ordre de $G$.

\subsection{Immersions ferm\'ees et quotients}

Nous commen\c cons par rappeler une propri\'et\'e (bien connue) de prolongement des $K$-sch\'emas en groupes finis.

\begin{lem}
\label{lemannexe}
Soit $G_K$ un $K$-sch\'ema en groupes fini, et soit $V\subseteq S$ un ouvert de $S$. Alors un prolongement de $G_K$ en un $V$-sch\'ema en groupes fini et lisse, s'il existe, est le mod\`ele de N\'eron de $G_K$ sur $V$ --- selon la terminologie de \cite[\S\/1.2, def. 1]{ray} --- donc est unique. En particulier, si l'ordre de $G_K$ est inversible sur $V$, alors il existe au plus un prolongement de $G_K$ en un $V$-sch\'ema en groupes fini et plat.
\end{lem}

\begin{lem}
\label{lem29}
Soit $H$ un $S$-sch\'ema en groupes fini et plat, et soit $f:H\rightarrow \A$ un morphisme de $S$-sch\'emas en groupes. Alors l'image de $f$ est un sous-groupe fini et plat de $\A$.
\end{lem}

\begin{proof}
L'image de $f$ est un sous-groupe quasi-fini et plat de $\A$. De plus, d'apr\`es \cite[corollaire 5.4.3 (ii)]{ega2}, l'image de $f$ est propre (car $\A$ est s\'epar\'e sur $S$, et $H$ est
propre sur $S$), donc l'image de $f$ est finie d'apr\`es le \og main theorem \fg{} de Zariski
(plus pr\'ecis\'ement, un morphisme quasi-fini et propre est fini \cite[cor. 4.4.11]{ega3}).
\end{proof}

\begin{dfn}
Soit $p$ un nombre premier. Nous dirons que $p$ est \emph{peu ramifi\'e} dans $S$ si, pour tout $\mathfrak{p}\in S$ de caract\'eristique r\'esiduelle $p$, l'indice de ramification absolu $e_{\mathfrak{p}}$ de $\mathfrak{p}$ satisfait l'in\'egalit\'e $e_{\mathfrak{p}}<p-1$.
\end{dfn}

\begin{rmq}
Si $p$ est inversible dans $S$, c'est-\`a-dire si $S$ ne contient aucun point de caract\'eristique r\'esiduelle $p$, alors $p$ est peu ramifi\'e dans $S$.
\end{rmq}

\begin{lem}
\label{lem210}
Avec les hypoth\`eses du lemme \ref{lem29}, supposons que $f_K:H_K\rightarrow A_K$ soit une immersion ferm\'ee, et que l'ordre de $H$ soit une puissance d'un nombre premier $p$ peu ramifi\'e dans $S$. Alors $f$ est une immersion ferm\'ee.
\end{lem}

\begin{proof}
Soit $H'$ l'image de $f$, et soit $V:=S[p^{-1}]\subseteq S$ le compl\'ementaire de l'ensemble des points de caract\'eristique $p$. Alors le morphisme naturel $H_V\rightarrow H'_V$ est un morphisme entre sch\'emas en groupes finis et lisses sur $V$, et sa fibre g\'en\'erique est un isomorphisme, donc c'est un isomorphisme, d'apr\`es le lemme \ref{lemannexe}. Soit d'autre part $\mathfrak{p}\in S$ un point ferm\'e de caract\'eristique $p$, et soit $S_{(\mathfrak{p})}=\Spec(R_{(\mathfrak{p})})$ o\`u $R_{(\mathfrak{p})}$ est le localis\'e de $R$ en $\mathfrak{p}$. Soit $e_{\mathfrak{p}}$ l'indice de ramification absolu de $\mathfrak{p}$, l'in\'egalit\'e $e_{\mathfrak{p}}<p-1$ implique que la fl\`eche $H\times_S S_{(\mathfrak{p})}\rightarrow H'\times_S S_{(\mathfrak{p})}$ est un isomorphisme, d'apr\`es \cite[\S\/7.5, lemma 5]{ray}. On en d\'eduit que $H$ et $H'$ sont isomorphes.
\end{proof}

\begin{ntt}
\label{notation1}
Nous noterons $\A^{\circ}$ la composante neutre de $\A$, et $\Phi:=\A/\A^{\circ}$ le groupe des composantes de $\A$. Si $\Gamma$ est un sous-groupe ouvert de $\Phi$, nous noterons $\mathcal{A}^{\Gamma}$ l'image r\'eciproque de $\Gamma$ par le morphisme canonique $\A\rightarrow\Phi$. Enfin, nous noterons $\A^t$ le mod\`ele de N\'eron de la $K$-vari\'ete ab\'elienne duale $A_K^t$ de $A_K$.

Nous fixons \`a pr\'esent un sous-groupe ouvert $\Gamma$ de $\Phi$, ainsi qu'un sous-groupe fini et plat $G$ de $\mathcal{A}^{\Gamma}$.
Nous noterons $B_K:=A_K/G_K$ la vari\'et\'e ab\'elienne quotient, $\B$ (resp. $\B^t$) le mod\`ele de N\'eron de $B_K$ (resp. $B_K^t$), et $\Psi$ le groupe des composantes de $\B$.
\end{ntt}

Par dualit\'e, nous avons une inclusion $G_K^D\subseteq B_K^t$. Dans certains cas, cette inclusion se prolonge en une immersion ferm\'ee $G^D\rightarrow \B^t$.

\begin{lem}
\label{lem21}
On suppose que $\A$ a bonne r\'eduction en toutes les places divisant $\ord(G)$. Alors l'inclusion $G_K^D\rightarrow B_K^t$ se prolonge en une immersion ferm\'ee $G^D\rightarrow \B^t$.
\end{lem}

\begin{proof}
Soit $U\subseteq S$ l'ouvert de bonne r\'eduction de $\A$, et soit
$$V:=S[\ord(G)^{-1}]\subseteq S$$
le compl\'ementaire de l'ensemble des points de caract\'eristique divisant $\ord(G)$. Comme $G^D$ est tu\'e par $\ord(G)$, on en d\'eduit que $G^D$ est \'etale (donc lisse) sur $V$, d'o\`u par la propri\'et\'e universelle du mod\`ele de N\'eron une fl\`eche
$$
\begin{CD}
G_V^D @>>> \B^t_V \\
\end{CD}
$$
qui prolonge l'inclusion $G_K^D\rightarrow B_K^t$, et qui est une immersion ferm\'ee d'apr\`es le lemme \ref{lemannexe}.

D'autre part, nous avons une suite exacte de faisceaux
$$
\begin{CD}
0 @>>> G_U @>>> \A_U @>>> \B_U @>>> 0 \\
\end{CD}
$$
pour la topologie fppf sur $U$, donc par dualit\'e des sch\'emas ab\'eliens un morphisme
$$
\begin{CD}
G_U^D @>>> \B^t_U \\
\end{CD}
$$
qui est une immersion ferm\'ee. Enfin, d'apr\`es l'hypoth\`ese ci-dessus, $U$ et $V$ forment un recouvrement ouvert de $S$. En recollant les morceaux, nous obtenons une fl\`eche
$$
\begin{CD}
G^D @>>> \B^t \\
\end{CD}
$$
qui est une immersion ferm\'ee.
\end{proof}

\begin{lem}
\label{lem22}
Le faisceau quotient $\A^{\Gamma}/G$ (pour la topologie fppf sur $S$) est repr\'esentable par un $S$-sch\'ema en groupes lisse, s\'epar\'e et de type fini, que nous noterons $B$.
\end{lem}

\begin{proof}
Comme $\A^{\Gamma}$ est de type fini sur $S$, et $S$ r\'egulier de dimension $1$, le faisceau quotient $\A^{\Gamma}/G$ est bien repr\'esentable d'apr\`es \cite[chap. IV, th\'eor\`eme 4.C]{anan}. La projection canonique $\varphi:\A^{\Gamma}\rightarrow B$ est fid\`element plate, et $\A^{\Gamma}$ est un $S$-sch\'ema
plat, donc $B$ est un $S$-sch\'ema plat. De plus $B$ est lisse sur
$S$ gr\^ace au crit\`ere de lissit\'e par fibres \cite[\S\/2.4, prop. 8]{ray}. On voit que $B$ est s\'epar\'e car $\A^{\Gamma}$ est s\'epar\'e et $G$ est ferm\'e dans $\A^{\Gamma}$.
\end{proof}

\begin{ntt}
\label{notation2}
Soit $\varphi_K:A_K\rightarrow B_K$ l'isog\'enie d\'efinissant $B_K$, et soit $\varphi^t_K:B_K^t\rightarrow A_K^t$ l'isog\'enie duale de $\varphi$. Nous noterons $\varphi$ et $\varphi^t$ les prolongements respectifs de ces isog\'enies sur les mod\`eles de N\'eron.
\end{ntt}

\begin{lem}
\label{lem24}
Supposons que les hypoth\`eses du lemme \ref{lem21} soient satisfaites, ou que $\A$ soit semi-stable. Alors $B$ est isomorphe \`a un sous-groupe ouvert de $\B$.
\end{lem}

\begin{proof}
Sous les hypoth\`eses envisag\'ees, l'ordre de $\varphi$ est premier aux caract\'eris\-tiques r\'esiduelles des places de mauvaise r\'eduction de $\A$, ou $\A$ est semi-stable. Par suite, le morphisme $\varphi:\A^{\circ}\rightarrow\B^{\circ}$ est surjectif, \`a noyau quasi-fini et plat, d'apr\`es \cite[\S\/7.3, prop. 6]{ray}. De m\^eme, le noyau de $\varphi:\A\rightarrow\B$ est un sous-groupe quasi-fini, plat et ferm\'e de $\A$, donc est \'egal \`a l'adh\'erence sch\'ematique dans $\A$ de sa fibre g\'en\'erique. Or cette derni\`ere n'est autre que $G$. On en d\'eduit le r\'esultat.
\end{proof}

\subsection{Accouplement de monodromie}
Nous rappelons ici quelques faits concernant l'accouplement (dit \og de monodromie \fg{}) introduit par Grothendieck dans \cite{gro7}.
Soit $\mathcal{P}_K$ le fibr\'e de Poincar\'e sur $A_K\times_K A_K^t$, c'est-\`a-dire le fibr\'e universel permettant d'identifier $A_K^t$ au foncteur de Picard de $A_K$. Gr\^ace au th\'eor\`eme du carr\'e, on peut munir $\mathcal{P}_K$
d'une unique structure de biextension de
$(A_K,A_K^t)$ par $\gmk$, que l'on appelle la biextension de
Weil, et que l'on note $W_K$.

Soit $\Phi'$ le groupe des composantes de $\A^t$. On
d\'eduit de la lecture de \cite[expos\'e VIII, th\'eor\`eme 7.1, b)]{gro7} un accouplement (associ\'e \`a $W_K$)
$$
\Phi\times_S\Phi'\longrightarrow ({\bf Q}/{\bf Z})_S
$$
qui repr\'esente l'obstruction \`a prolonger $W_K$ en une biextension de $(\A,\A^t)$ par
$\gm$. Plus pr\'ecis\'ement, nous avons la proposition suivante.

\begin{prop}
\label{woz}
Soit $M$ (resp. $M'$) un sous-groupe de $\Phi$ (resp. $\Phi'$). Alors il
existe une
(unique)
biextension $W$ de $(\A^{M},\A^{t,M'})$ par
$\gm$ prolongeant la biextension de Weil $W_K$ sur
$(A_K,A_K^t)$ si et
seulement
si $M$ et $M'$ sont orthogonaux sous
l'accouplement.
\end{prop}

\begin{proof}
Le r\'esultat d\'ecoule de \cite[expos\'e VIII, th\'eor\`eme 7.1, b)]{gro7} dans
le cas particulier o\`u la base est un trait. En outre, apr\`es
lecture de \cite[expos\'e VIII, remarque 7.2]{gro7}, on voit que ce
r\'esultat s'\'etend au spectre d'un anneau de Dedekind, ce qui est bien le cas ici.
\end{proof}

\begin{ntt}
\label{notation3}
On sait que l'image de $\A^{\circ}$ par $\varphi:\A\rightarrow\B$ est \`a valeurs dans $\B^{\circ}$. Par passage au quotient, on en d\'eduit une application $\overline\varphi:\Phi\rightarrow\Psi$ entre les groupes de composantes. Dans la suite, nous noterons
$$
\Lambda:=\overline\varphi(\Gamma)
$$
et $\Lambda'$ d\'esignera l'orthogonal de $\Lambda$ sous l'accouplement de monodromie (associ\'e \`a la biextension de Weil $W'_K$ de $(B_K,B_K^t)$ par $\gm$).
\end{ntt}

\begin{lem}
\label{lem26}
Avec les notations pr\'ec\'edentes, $G$ est un sous-groupe de $\ker \varphi$, la fl\`eche $B\rightarrow \B$ est \`a valeurs dans $\B^{\Lambda}$, $\ker \varphi^t$ est contenu dans $\B^{t,\Lambda'}$, et $\overline\varphi^t(\Lambda')\subseteq \Gamma'$.
\end{lem}

\begin{proof}
Il est clair que $G$ et $\ker\varphi$ ont la m\^eme fibre g\'en\'erique. De plus, ces deux groupes sont ferm\'es dans $\A$, et $G$ est \'egal \`a l'adh\'erence sch\'ematique dans $\A$ de sa fibre g\'en\'erique (ceci est d\^u au fait que $G$ est plat sur $S$), d'o\`u l'inclusion de $G$ dans $\ker \varphi$. Gr\^ace \`a la propri\'et\'e universelle du quotient, on en d\'eduit un morphisme $B\rightarrow \B^{\Lambda}$.

Pour obtenir l'inclusion $\ker \varphi^t\subseteq\B^{t,\Lambda'}$, il suffit de montrer que $\ker \overline\varphi^t\subseteq\Lambda'$.
Soient respectivement $\Phi'$ et $\Psi'$ les groupes de composantes de $\A^t$ et $\B^t$. On peut r\'esumer la situation par un diagramme
\[
\begin{CD}
\Phi @.\times_S @. \Phi' @>>> ({\bf Q}/{\bf Z})_S \\
@V\overline\varphi VV @. @AA\overline\varphi^t A @| \\
\Psi @.\times_S @. \Psi' @>>> ({\bf Q}/{\bf Z})_S \\
\end{CD}
\]
dans lequel on note $<,>_{\A}$ l'accouplement du haut (associ\'e \`a
$W_K$), et
$<,>_{\B}$
l'accouplement du bas (associ\'e \`a
$W_K'$). Il r\'esulte alors de \cite[expos\'e VIII, 7.3.1]{gro7} et de
l'identit\'e
$(\varphi_K\times {\rm
id}_{B_K^t})^*(W_K')=({\rm id}_{A_K}\times
\varphi_K^t)^*(W_K)$
que le
diagramme pr\'ec\'edent est commutatif,
c'est-\`a-dire que nous avons, pour tout $x\in \Phi$ et tout $y\in \Psi'$,
l'\'egalit\'e
$$<\overline\varphi(x),y>_{\B}=<x,\overline\varphi^t(y)>_{\A}\,.$$
Cette identit\'e permet aussit\^ot alors de voir que
l'orthogonal de $\Lambda:=\overline\varphi(\Gamma)$ sous l'accouplement $<,>_{\B}$
contient le noyau de $\overline\varphi^t$. L'inclusion $\overline\varphi^t(\Lambda')\subseteq \Gamma'$ s'obtient \'egalement en observant la m\^eme identit\'e.
\end{proof}


\section{Homomorphisme de classes et restrictions}

\subsection{Suites exactes et dualit\'e}

Le \og petit site fppf \fg{} sur $S$ est la cat\'egorie des sch\'emas plats sur $S$
munie d'une structure de site pour la topologie fppf. Nous appellerons
faisceau (pour la topologie fppf) sur $S$ un faisceau sur ce site.
Nous allons travailler \`a pr\'esent dans la cat\'egorie des faisceaux ab\'eliens sur le petit site fppf de $S$.

\begin{lem}
Soit $H$ un $S$-sch\'ema en groupes fini et plat. On se donne un morphisme $f:H\rightarrow \A$ dont la fibre g\'en\'erique $f_K:H_K\rightarrow A_K$ est une immersion ferm\'ee. Alors $f$ est un monomorphisme de faisceaux sur le petit site fppf de $S$.
\end{lem}

\begin{proof}
Soit en effet $T\rightarrow S$ un $S$-sch\'ema plat. Nous avons un diagramme commutatif
$$
\begin{CD}
H(T) @>>> \A(T) \\
@VVV @VVV \\
H_K(T_K) @>>> A_K(T_K) \\
\end{CD}
$$
dans lequel la fl\`eche verticale de gauche est injective, $H$ \'etant s\'epar\'e sur $S$, et $T_K$ \'etant dense dans $T$ (gr\^ace \`a la platitude de $T$). De plus, la fl\`eche horizontale du bas est injective puisque $f_K$ est une immersion. Donc la fl\`eche horizontale du haut est injective.
\end{proof}

Nous sommes \`a pr\'esent en mesure de donner une construction de l'homomorphisme de classes $\psi$. Cette construction diff\`ere l\'eg\`erement de celle donn\'ee par l'auteur dans \cite{gil1} et \cite{gil2}. Elle est cependant bas\'ee sur les m\^emes principes.

\begin{lem}
\label{central}
Avec les notations \ref{notation1}, \ref{notation2} et \ref{notation3}, supposons que l'on ait simultan\'ement une immersion ferm\'ee $G\rightarrow \A^{\Gamma}$ et un morphisme $G^D\rightarrow \B^{t,\Lambda'}$ qui prolonge l'immersion canonique $G_K^D\rightarrow B_K^t$.

Soient $B:=\A^{\Gamma}/G$ et $F:=\B^{t,\Lambda'}/G^D$ les faisceaux quotients respectifs (pour la topologie fppf sur $S$).
Alors nous avons un diagramme commutatif
\begin{equation}
\label{deuxsuites}
\begin{CD}
0 @>>> G @>>> \A^{\Gamma} @>>> B @>>> 0 \\
@. @| @V\alpha VV @V\beta VV \\
0 @>>> G @>>> \Ext(F,\gm) @>>> \Ext(\B^{t,\Lambda'},\gm) @>>> 0 \\
\end{CD}
\end{equation}
dans lequel les fl\`eches $\alpha$ est $\beta$ sont d\'efinies en termes de biextensions, et la suite exacte du bas est obtenue en appliquant le foncteur $\Hom(-,\gm)$ \`a la suite exacte
\begin{equation}
\label{suitebis}
\begin{CD}
0 @>>> G^D @>>> \B^{t,\Lambda'} @>>> F @>>> 0 \\
\end{CD}
\end{equation}
\end{lem}

\begin{proof}
D'apr\`es \cite[lemme 2.4]{gil1}, le faisceau $\Hom(\B^{t,\Lambda'},\gm)$ est nul. De plus, il d\'ecoule de \cite[expos\'e VIII, 3.3.1]{gro7} que le faisceau $\Ext(G,\gm)$ est nul. Enfin le faisceau $\Hom(G^D,\gm)$ est repr\'esentable par $G$. Ainsi la suite exacte du bas dans le diagramme (\ref{deuxsuites}) d\'ecoule de la suite exacte (\ref{suitebis}) par application du foncteur $\Hom(-,\gm)$.

Soit $H$ l'image du morphisme $G^D\rightarrow\B^{t,\Lambda'}$, et soit $A^t$ le quotient $\B^{t,\Lambda'}/H$, lequel est repr\'esentable par un $S$-sch\'ema en groupes lisse (lemme \ref{lem22}), contrairement \`a $F$ qui en g\'en\'eral ne l'est pas. Alors on dispose d'un morphisme canonique de faisceaux $F\rightarrow A^t$. En outre, d'apr\`es le lemme \ref{lem26}, on dispose d'un morphisme canonique $A^t\rightarrow \A^{t,\overline\varphi^t(\Lambda')}$, et, toujours par le lemme \ref{lem26}, $\overline\varphi^t(\Lambda')$ est un sous-groupe de $\Gamma'$. En recollant tous les morceaux, on obtient un morphisme de faisceaux $F\rightarrow \A^{t,\Gamma'}$.

Soit $W$ (resp. $W'$) la biextension de $(\A^{\Gamma},\A^{t,\Gamma'})$ (resp. de $(\B^{\Lambda},\B^{t,\Lambda'})$) par $\gm$ qui prolonge la biextension de Weil.  On d\'efinit alors $\alpha$ comme \'etant le compos\'e des morphismes
$$
\begin{CD}
\alpha:\A^{\Gamma} @>>> \Ext(\A^{t,\Gamma'},\gm) @>>> \Ext(F,\gm) \\
\end{CD}
$$
o\`u la premi\`ere fl\`eche est d\'efinie par la biextension $W$. De m\^eme, on dispose d'un morphisme canonique $B\rightarrow\B^{\Lambda}$, et on d\'efinit $\beta$ comme \'etant le compos\'e des morphismes
$$
\begin{CD}
\beta:B @>>> \B^{\Lambda} @>>> \Ext(\B^{t,\Lambda'},\gm) \\
\end{CD}
$$
et il est clair que le diagramme (\ref{deuxsuites}) commute.
\end{proof}

\begin{prop}
\label{classinvariant}
Avec les notations \ref{notation1}, \ref{notation2} et \ref{notation3}, supposons que l'on ait simultan\'ement une immersion ferm\'ee $G\rightarrow \A$, et un morphisme $i:G^D\rightarrow \B^t$ qui prolonge l'immersion canonique $G_K^D\rightarrow B_K^t$. Alors les hypoth\`eses du lemme \ref{central} sont satisfaites. Supposons de plus que le cobord $\Delta$ associ\'e \`a la suite exacte
$$
\begin{CD}
0 @>>> G @>>> \A @>>> B @>>> 0 \\
\end{CD}
$$
soit surjectif. Alors le morphisme
$$
\begin{CD}
\B^{\Lambda}(S) @>\gamma >> \ext(\B^{t,\Lambda'},\gm) @>\delta >> H^1(S,G) \\
\end{CD}
$$
est surjectif, $\delta$ d\'esignant le cobord associ\'e \`a la suite du bas dans le diagramme $(\ref{deuxsuites})$, et $\gamma$ \'etant d\'efini gr\^ace \`a la biextension $W'$. Enfin, l'homomorphisme de classes $\psibemol$ d\'efini comme \'etant la compos\'ee des fl\`eches suivantes
$$
\begin{CD}
\B^{\Lambda}(S) @>\delta\circ\gamma >> H^1(S,G) @>\pi >> \rpic(G^D) \\
\end{CD}
$$
est \'egalement surjectif. En termes de fibr\'es en droites, ceci \'equivaut \`a la surjectivit\'e de l'application
\begin{equation}
\begin{split}
\B^{\Lambda}(S)\;\longrightarrow & \;\rpic(G^D)\\
(x:S\rightarrow \B^{\Lambda})\;\longmapsto & \;(x\times i)^*(t(W'))\\
\end{split}
\end{equation}
o\`u $t(W')$ d\'esigne le $\gm$-torseur sur $\B^{\Lambda}\times_S\B^{t,\Lambda'}$ associ\'e \`a la biextension $W'$.
\end{prop}

\begin{proof}
Le lemme \ref{lem26} montre que les hypoth\`eses du lemme \ref{central} sont satisfaites, en prenant $\Gamma=\Phi$.
Nous avons un diagramme commutatif
$$
\begin{CD}
B(S) @>\Delta >> H^1(S,G) \\
@V\beta VV @| \\
\ext(\B^{t,\Lambda'},\gm) @>\delta >> H^1(S,G) \\
\end{CD}
$$
entre les deux cobords des suites exactes du diagramme (\ref{deuxsuites}). D'autre part, $\beta$ se d\'efinit comme \'etant la compos\'ee
$$
\begin{CD}
\beta:B(S) @>>> \B^{\Lambda}(S) @>\gamma >> \ext(\B^{t,\Lambda'},\gm) \\
\end{CD}
$$
o\`u $\gamma$ est d\'efinie \`a l'aide de la biextension $W'$. Par suite, la surjectivit\'e de $\Delta$ entra\^ine celle de $\delta\circ\beta$ ainsi que celle de $\delta\circ\gamma$, et \emph{a fortiori} celle de $\psibemol$. Le dernier point d\'ecoule de la description (dite \og description g\'eom\'etrique \fg{}) de l'homomorphisme de classes $\psibemol$ en termes de restriction du fibr\'e $t(W')$. Nous renvoyons pour cela le lecteur aux travaux pr\'ec\'edents de l'auteur (voir \cite[lemme 3.2]{gil1}), qui s'adaptent sans peine \`a la situation que nous consid\'erons ici.
\end{proof}


\subsection{Restriction de Weil}

Soit $K'$ une extension finie de $K$, et soit $S':=\Spec(R')$ o\`u $R'$ est la cl\^oture int\'egrale de $R$ dans $K'$.
Soit $X_{K'}$ une $K'$-vari\'et\'e ab\'elienne. Nous noterons $\X_{/S'}\rightarrow S'$ le mod\`ele de N\'eron de $X_{K'}$ sur $S'$.

Rappelons la proposition suivante (voir \cite[Prop. 4.1]{edix}).

\begin{prop}
La restriction de Weil $\Re_{K'/K}(X_{K'})$ est repr\'esentable par une $K$-vari\'et\'e ab\'elienne, et la restriction de Weil $\Re_{S'/S}(\X_{/S'})$ est repr\'esentable par le mod\`ele de N\'eron de la vari\'et\'e ab\'elienne $\Re_{K'/K}(X_{K'})$.
\end{prop}

\begin{rmq}
Le r\'esultat ci-dessous ne sera pas utilis\'e dans nos constructions. Il permet cependant de se faire une id\'ee des propri\'et\'es de r\'eduction des vari\'et\'es ab\'eliennes que nous allons manipuler dans la suite.
\end{rmq}

\begin{lem}
\label{lem16}
Soit $\mathfrak{p}\in S$ un point ferm\'e de $S$ qui ne se ramifie pas dans $K'$, et soit $\mathfrak{q}$ un point de $S'$ au-dessus de $\mathfrak{p}$. Alors $\Re_{S'/S}(\X_{/S'})$ a bonne (resp. semi-stable) r\'eduction en $\mathfrak{p}$ si et seulement si $\X_{/S'}$ a bonne (resp. semi-stable) r\'eduction en $\mathfrak{q}$.
\end{lem}

\begin{proof}
Soit $k(\mathfrak{p})$ le corps r\'esiduel de $\mathfrak{p}$. Soient $\mathfrak{q}_1,\dots,\mathfrak{q}_r$ les point de $S'$ au-dessus de $\mathfrak{p}$, et soient $k(\mathfrak{q}_1),\dots,k(\mathfrak{q}_r)$ leurs corps r\'esiduel respectifs. Soit $T:=S'\times_S k(\mathfrak{p})$, alors, comme $\mathfrak{p}$ ne se ramifie pas dans $K'$, nous avons $T=\Spec(k(\mathfrak{q}_1)\times\dots\times k(\mathfrak{q}_r))$. Nous avons d'une part, par compatibilit\'e de la restriction de Weil avec le changement de base,
$$
\Re_{S'/S}(X_{/S'})\times_S k(\mathfrak{p})=\Re_{T/k(\mathfrak{p})}(X_{/S'}\times_{S'} T)
$$
et d'autre part, sachant que $T=\amalg_{i=1}^r\Spec(k(\mathfrak{q}_i))$, il est facile de voir que
$$
\Re_{T/k(\mathfrak{p})}(X_{/S'}\times_{S'} T)=\Pi_{i=1}^r\Re_{k(\mathfrak{q}_i)/k(\mathfrak{p})}(X_{/S'}\times_{S'} k(\mathfrak{q}_i)).
$$
On se sert alors du fait que la restriction d'une vari\'et\'e ab\'elienne (resp. semi-ab\'elienne) sur un corps est une vari\'et\'e ab\'elienne (resp. semi-ab\'elienne).
\end{proof}


\subsection{Sous-groupes finis des restrictions}

Dor\'enavant, $G$ d\'esigne un $S$-sch\'ema en groupes fini et plat.

\begin{lem}
\label{lem32}
Soit $A_K:=\Re_{K'/K}(X_{K'})$. Supposons que l'on dispose d'un morphisme
$$
\begin{CD}
G_{S'} @>>> \X_{/S'} \\
\end{CD}
$$
qui soit une inclusion sur les fibres g\'en\'eriques. Alors en composant les fl\`eches
$$
\begin{CD}
G @>>> \Re_{S'/S}(G_{S'}) @>>> \Re_{S'/S}(\X_{/S'})=\A \\
\end{CD}
$$
on obtient un morphisme $G\rightarrow \A$. Si en outre $G_{S'}\rightarrow \X_{/S'}$ est une immersion ferm\'ee, ou si l'ordre de $G$ est une puissance d'un nombre premier peu ramifi\'e dans $S$, alors $G\rightarrow \A$ est une immersion ferm\'ee.
\end{lem}

\begin{proof}
La d\'efinition du morphisme $G\rightarrow \A$ est claire. Si $G_{S'}\rightarrow \X_{/S'}$ est une immersion ferm\'ee, alors $\Re_{S'/S}(G_{S'})\rightarrow \Re_{S'/S}(\X_{/S'})$ est une immersion ferm\'ee d'apr\`es \cite[\S\/7.6, prop. 2]{ray}. D'autre part, $G\rightarrow\Re_{S'/S}(G_{S'})$ est une immersion ferm\'ee d'apr\`es \cite[\S\/7.6, p. 197]{ray}, le groupe $G$ \'etant affine, donc s\'epar\'e. Ceci montre le premier point. Le second point d\'ecoule du lemme \ref{lem210}.
\end{proof}

\begin{lem}
\label{lem46}
Supposons que $X_{K'}$ provienne d'une $K$-vari\'et\'e ab\'elienne $X_K$, et que $G\subseteq\X$ soit un sous-groupe du mod\`ele de N\'eron $\X$ de $X_K$ sur $S$. Supposons que $\X$ ait r\'eduction semi-stable en les points de $S$ de caract\'eristique r\'esiduelle divisant $\ord(G)$. Alors le morphisme compos\'e
$$
\begin{CD}
G_{S'} \subseteq \X_{S'} @>>> \X_{/S'} \\
\end{CD}
$$
est une immersion ferm\'ee. Avec les notations du lemme \ref{lem32}, le morphisme $G\rightarrow \A$ est alors une immersion ferm\'ee.
\end{lem}

\begin{proof}
Soit $V':=S'[\ord(G)^{-1}]\subseteq S'$ le compl\'ementaire de l'ensemble des points de caract\'eristique divisant $\ord(G)$, de sorte que $G_{V'}$ est \'etale sur $V'$. Par suite, nous avons une immersion ferm\'ee $G_{V'}\rightarrow \X_{/V'}$.
Soit d'autre part $U'$ le plus grand ouvert de $S'$ sur lequel $\X_{S'}$ a r\'eduction semi-stable. Alors la fl\`eche $\X_{U'}\rightarrow \X_{/U'}$ est une immersion ouverte d'apr\`es \cite[\S\/7.4, prop. 3]{ray}. Par suite, la fl\`eche $G_{U'}\rightarrow \X_{/U'}$ qui s'en d\'eduit est une immersion, comme compos\'ee de deux immersions. Comme, par hypoth\`ese, $U'$ et $V'$ constituent un recouvrement ouvert de $S'$, on en d\'eduit le r\'esultat. Le dernier point d\'ecoule imm\'ediatement du lemme \ref{lem32}.
\end{proof}


\subsection{L'homomorphisme de classes surjectif}
\label{surjectif}

Soit $X_{K'}$ une $K'$-vari\'et\'e ab\'elienne, $\X_{/S'}$ son mod\`ele de N\'eron sur $S'$, et $\A=\Re_{S'/S}(\X_{/S'})$. On se donne un morphisme $G_{S'}\rightarrow \X_{/S'}$ tel que le morphisme induit $G\rightarrow \A$ soit une immersion ferm\'ee.

\begin{thm}
\label{mainthm}
Supposons que la fl\`eche naturelle de changement de base
$$
\begin{CD}
H^1(S,G) @>>> H^1(S',G_{S'}) \\
\end{CD}
$$
soit nulle. Alors le morphisme cobord
$$
\begin{CD}
B(S) @>\Delta >> H^1(S,G)
\end{CD}
$$
est surjectif.
\end{thm}

\begin{proof}
Le point de d\'epart est la suite exacte
$$
\begin{CD}
0 @>>> G @>>> \A @>>> B @>>> 0\\
\end{CD}
$$
d'o\`u l'on d\'eduit une suite exacte
$$
\begin{CD}
\A(S) @>>> B(S) @>\Delta>> H^1(S,G) @>>> H^1(S,\A). \\
\end{CD}
$$

Nous avons d'autre part un isomorphisme canonique
$$
H^1(S,\A)=H^1(S,\Re_{S'/S}(\X_{/S'}))\simeq H^1(S',\X_{/S'}).
$$
En effet, les $S$-groupes consid\'er\'es \'etant lisses, on peut consid\'erer les groupes de cohomologie comme \'etant calcul\'es dans la topologie \'etale. On se sert alors du fait que $R^1h_*=0$ en cohomologie \'etale, $h:S'\rightarrow S$ \'etant le morphisme (fini) de changement de base.

Pour montrer que $\Delta$ est surjectif, il suffit de montrer que la fl\`eche
$$
\begin{CD}
q:H^1(S,G) @>>> H^1(S',\X_{/S'}) \\
\end{CD}
$$
est nulle. Or cette fl\`eche $q$ est en fait obtenue en composant les fl\`eches
$$
\begin{CD}
H^1(S,G) @>>> H^1(S',G_{S'})@>>> H^1(S',\X_{/S'}) \\
\end{CD}
$$
o\`u la premi\`ere fl\`eche est obtenue par changement de base, et la seconde est induite par le morphisme $G_{S'}\rightarrow \X_{/S'}$. On en d\'eduit le r\'esultat.
\end{proof}

\begin{lem}
\label{nulchbase}
Il existe une extension $K'$ de $K$, non ramifi\'ee en dehors des places divisant $\ord(G)$, telle que la fl\`eche
$$
\begin{CD}
H^1(S,G) @>>> H^1(S',G_{S'}) \\
\end{CD}
$$
soit nulle.
\end{lem}

\begin{proof}
Rappelons un fait bien connu : le groupe $H^1(S,G)$ est un groupe ab\'elien fini. Donc l'image de l'application (injective) de restriction \`a la fibre g\'en\'erique
$$
\begin{CD}
H^1(S,G) @>>> H^1(K,G_K) \\
\end{CD}
$$
est finie. De plus, chaque \'el\'ement de cette image se trivialise dans une extension de degr\'e fini de $K$, non ramifi\'ee en dehors des places divisant $\ord(G)$ (en effet, si l'on inverse les places divisant $\ord(G)$, alors le groupe $G$, ainsi que ses torseurs, deviennent \'etales). Soit $K'$ le compositum de ces extensions. Alors $K'$ est une extension finie de $K$, non ramifi\'ee en dehors des places divisant $\ord(G)$, et la fl\`eche
$$
\begin{CD}
H^1(S,G) @>>> H^1(S',G_{S'}) \\
\end{CD}
$$
est nulle, ce qu'on voulait.
\end{proof}

\begin{proof}[D\'emonstration du th\'eor\`eme \ref{sgsemistable}]
Soit $K'$ une extension de $K$ telle que la fl\`eche de changement de base $H^1(S,G)\rightarrow H^1(S',G_{S'})$ soit nulle (voir le lemme \ref{nulchbase}). En appliquant le lemme \ref{lem46}, nous obtenons une immersion ferm\'ee $G\rightarrow \A=\Re_{S'/S}(\X_{/S'})$. Les hypoth\`eses du th\'eor\`eme \ref{mainthm} sont donc satisfaites, et le r\'esultat en d\'ecoule.
\end{proof}


\section{Cons\'equences et exemples}

\subsection{Le cas o\`u $G^D$ est un groupe constant cyclique}

Nous fixons dans la suite un nombre premier $p$, et un entier naturel $n$.

\begin{prop}
\label{prop41}
Supposons qu'il existe une immersion ferm\'ee $\mu_{p^n}\rightarrow\X$, o\`u $\X$ est le mod\`ele de N\'eron d'une $K$-vari\'et\'e ab\'elienne $X_K$, ayant r\'eduction semi-stable en les points de $S$ de caract\'eristique $p$. Alors il existe une extension $K'$ de $K$ telle que les hypoth\`eses de la proposition \ref{classinvariant} soient satisfaites pour $\mu_{p^n}\rightarrow\A=\Re_{S'/S}(\X_{/S'})$.
\end{prop}

\begin{proof}
Soit $K'$ une extension de $K$ telle que la fl\`eche
$$
\begin{CD}
H^1(S,\mu_{p^n}) @>>> H^1(S',\mu_{p^n,S'}) \\
\end{CD}
$$
soit nulle (une telle extension existe d'apr\`es le lemme \ref{nulchbase}). D'apr\`es le lemme \ref{lem46}, le morphisme $\mu_{p^n}\rightarrow\A$ qui s'en d\'eduit est une immersion ferm\'ee. De plus, comme $\zpnz_S$ est \'etale sur $S$, on dispose d'un morphisme $\zpnz_S\rightarrow \B^t$ qui prolonge l'immersion canonique sur les fibres g\'en\'eriques. Enfin le th\'eor\`eme \ref{mainthm} s'applique, et implique que toutes les hypoth\`eses de la proposition \ref{classinvariant} sont satisfaites.
\end{proof}

Le corollaire qui suit est motiv\'e entre autres par la question pos\'ee par Agboola et Pappas dans \cite{ap0}, m\^eme s'il ne fournit pas de r\'eponse \`a cette derni\`ere. On peut se demander quelles classes dans $\pic(S)$ peuvent s'obtenir comme le pull-back (par un $S$-point) d'un fibr\'e en droites sur une vari\'et\'e arithm\'etique donn\'ee. Ici, nous montrons que tout \'el\'ement de $\pic(S)[p^n]$ est le pull-back d'un fibr\'e en droites rigidifi\'e sur un sous-groupe ouvert du mod\`ele de N\'eron d'une vari\'et\'e ab\'elienne convenable.

\begin{cor}
\label{petitcor}
Il existe une $K$-vari\'et\'e ab\'elienne $B_K$, un sous-groupe $\Lambda$ du groupe des composantes du mod\`ele de N\'eron $\B$ de $B_K$, et un point de $p^n$-torsion $y\in \B^{t,\Lambda'}(S)$, tels que l'application
\begin{equation*}
\begin{split}
\B^{\Lambda}(S)\;\longrightarrow & \;\pic(S)[p^n]\\
(x:S\rightarrow \B^{\Lambda})\;\longmapsto & \;(x\times y)^*(t(W'))\\
\end{split}
\end{equation*}
soit surjective (avec les notations de la prop. \ref{classinvariant}).
\end{cor}

\begin{proof}
Il existe une $\mathbb{Q}$-vari\'et\'e ab\'elienne $C$, ayant bonne r\'eduction au-dessus de $p$, et dont le mod\`ele de N\'eron poss\`ede un sous-groupe isomorphe \`a $\mu_{p^n}$ (voir \cite[Th\'eor\`eme 3.1.1]{bbm}). Si $\C$ d\'esigne le mod\`ele de N\'eron de $C_K$ sur $S$ on obtient, gr\^ace au lemme \ref{lem46}, une immersion ferm\'ee $\mu_{p^n}\rightarrow \C$. On conclut en appliquant les propositions \ref{prop41} et \ref{classinvariant} et en se servant de l'isomorphisme $\rpic(\zpnz_S)\simeq\pic(S)[p^n]$.
\end{proof}

On en d\'eduit le th\'eor\`eme \ref{untrois} de la fa\c con suivante : on d\'ecompose $N$ en produit de puissances de nombres premiers distincts, puis on applique le corollaire \ref{petitcor} \`a chacune de ces puissances et l'on consid\`ere le produit des vari\'et\'es ainsi obtenues.

\begin{thm}
\label{essentiel}
Supposons que $p$ soit peu ramifi\'e dans $S$. Soit $K'$ une extension de $K$ telle que la fl\`eche de changement de base $H^1(S,\mu_p)\rightarrow H^1(S',\mu_{p,S'})$ soit nulle, et soit $E_{K'}$ une $K'$-courbe elliptique, ayant bonne r\'eduction au-dessus de $p$, et telle que $E_{K'}(K')$ poss\`ede un point d'ordre $p$. Soit $X_{K'}:=E_{K'}/\zpz_{K'}$, et soit $\A:=\Re_{S'/S}(\X_{/S'})$. Alors on dispose d'une immersion ferm\'ee $\mu_p\rightarrow \A$ satisfaisant les hypoth\`eses de la proposition \ref{classinvariant}. Enfin, le rang du groupe $B(S)$ est \'egal au rang du groupe $E_{K'}(K')$.
\end{thm}

\begin{proof}
Par la propri\'et\'e universelle du mod\`ele de N\'eron le point d'ordre $p$ dans $E_{K'}(K')$ d\'efinit un morphisme $\zpz_{S'}\rightarrow \E_{/S'}$ dont nous noterons $H$ l'image. D'apr\`es le lemme \ref{lem29}, $H$ est un sous-groupe fini et plat de $\E_{/S'}$. Soit $X_{K'}:=(E_{K'}/\zpz_{K'})^t$, alors, comme $E_{K'}$ a bonne r\'eduction au-dessus de $p$, le lemme \ref{lem21} s'applique et $H^D$ est un sous-groupe de $\X_{/S'}$. De plus, nous avons un morphisme $H^D\rightarrow \mu_{p,S'}$ d\'eduit du morphisme $\zpz_{S'}\rightarrow H$. On en d\'eduit un morphisme
$$
\begin{CD}
\nu:\Re_{S'/S}(H^D) @>>> \Re_{S'/S}(\mu_{p,S'})
\end{CD}
$$
et je pr\'etends que c'est un isomorphisme. En vertu de l'in\'egalit\'e $e<p-1$, il suffit de montrer que sa restriction \`a $V:=S[p^{-1}]$ est un isomorphisme. Soit $V'=S'[p^{-1}]$ alors la fl\`eche $\zpz_{V'}\rightarrow H_{V'}$ est un isomorphisme, donc $H^D_{V'}\rightarrow \mu_{p,V'}$ est un isomorphisme et par restriction la fl\`eche
$$
\begin{CD}
\Re_{V'/V}(H^D_{V'}) @>>> \Re_{V'/V}(\mu_{p,V'})
\end{CD}
$$
est un isomorphisme. Or nous avons $\Re_{V'/V}(Y_{V'})=\Re_{S'/S}(Y) \times_S V$ pour tout $S'$-sch\'ema $Y$, d'o\`u le r\'esultat.

L'isomorphisme $\nu$ donne naissance \`a une immersion ferm\'ee
$$
\begin{CD}
\mu_p \subseteq \Re_{S'/S}(\mu_{p,S'}) @>\nu^{-1}>> \Re_{S'/S}(H^D) @>>> \Re_{S'/S}(\X_{/S'})=\A
\end{CD}
$$
qui satisfait, comme nous allons le voir, les hypoth\`eses de la proposition \ref{classinvariant}. D'abord, le groupe $\zpz$ \'etant lisse sur $S$, il existe un morphisme $\zpz\rightarrow \B^t$ qui prolonge l'inclusion sur les fibres g\'en\'eriques. Ce morphisme est d'ailleurs une immersion ferm\'ee d'apr\`es le lemme \ref{lem210}.

Montrons \`a pr\'esent la surjectivit\'e du cobord. Le principe est le m\^eme que dans la preuve du th\'eor\`eme \ref{mainthm}, il suffit en fait de montrer que la fl\`eche
$$
\begin{CD}
q:H^1(S,\mu_p) @>>> H^1(S',\X_{S'}) \\
\end{CD}
$$
est nulle. Or cette fl\`eche $q$ est en fait obtenue en composant les fl\`eches
$$
\begin{CD}
H^1(S,\mu_p) @>>> H^1(S,\Re_{S'/S}(H^D)) @>>> H^1(S',H^D) @>>> H^1(S',\X_{S'}). \\
\end{CD}
$$
Nous avons d'autre part un diagramme commutatif
$$
\begin{CD}
H^1(S,\mu_p) @>>> H^1(S,\Re_{S'/S}(H^D)) \\
@VVV @VVV \\
H^1(S',\mu_{p,S'}) @<<< H^1(S',H^D) \\
\end{CD}
$$
dans lequel la fl\`eche verticale de gauche est nulle, et la fl\`eche horizontale du bas est injective, d'apr\`es le lemme \ref{sorite} ci-dessous. On en d\'eduit que $q$ est nulle.

Enfin, au sujet du rang, nous avons plus pr\'ecis\'ement les \'egalit\'es suivantes
$$
\rk(B(S))=\rk(\A(S))=\rk(\X_{/S'}(S'))=\rk(X_{K'}(K'))=\rk(E_{K'}(K')).
$$
La premi\`ere d\'ecoule de la finitude des groupes $\mu_p(S)$ et $H^1(S,\mu_p)$. Les deux suivantes d\'ecoulent des \'egalit\'es $\A(S)=\X_{/S'}(S')$ (vraie par d\'efinition m\^eme de la restriction de Weil) et $\X_{/S'}(S')=\X_K(K')$ (vraie par la propri\'et\'e de N\'eron). La derni\`ere est due au fait que $E_{K'}$ et $X_{K'}$ sont isog\`enes.
\end{proof}

\begin{lem}
\label{sorite}
Soit $M\rightarrow N$ un morphisme entre $S$-sch\'emas en groupes finis et plats, qui induise un isomorphisme sur les fibres g\'en\'eriques.
Alors le morphisme induit
$$
\begin{CD}
H^1(S,M) @>>> H^1(S,N) \\
\end{CD}
$$
est injectif.
\end{lem}

\begin{proof}
Il suffit d'examiner le diagramme commutatif
$$
\begin{CD}
H^1(S,M) @>>> H^1(K,M_K) \\
@VVV @VVV \\
H^1(S,N) @>>> H^1(K,N_K) \\
\end{CD}
$$
dans lequel les fl\`eches horizontales sont obtenues par restriction \`a la fibre g\'en\'erique, donc sont injectives. De plus, la fl\`eche verticale de droite est un isomorphisme (car $M_K\rightarrow N_K$ est un isomorphisme). On en d\'eduit que la fl\`eche verticale de gauche est injective.
\end{proof}


\subsection{Le cas quadratique imaginaire}

Dans cette section, on suppose que $R=\mathcal{O}_K$ est l'anneau des entiers de $K$. Dans le cas o\`u $K$ est un corps quadratique imaginaire et o\`u $p\geq 3$, trivialiser les $\mu_p$-torseurs revient \`a trivialiser les \'el\'ements de $\pic(\mathcal{O}_K)[p]$.

\begin{lem}
\label{naif}
Soit $p\geq 3$ un nombre premier, et soit $K\neq\mathbb{Q}(\sqrt{-3})$ un corps quadratique imaginaire. Soit $K'$ une extension de $K$ dans laquelle les \'el\'ements de $\pic(\mathcal{O}_K)[p]$ capitulent. Alors la fl\`eche
$$
\begin{CD}
H^1(S,\mu_p) @>>> H^1(S',\mu_{p,S'}) \\
\end{CD}
$$
est nulle.
\end{lem}

\begin{proof}
Comme $K$ est quadratique imaginaire, le groupe des unit\'es $\mathcal{O}_K^{\times}$ est un groupe fini. De plus, sous les hypoth\`eses envisag\'ees, $\mathcal{O}_K^{\times}$ est sans $p$-torsion : si $p>3$ c'est clair, et si $p=3$ le seul corps quadratique contenant les racines $3$-i\`emes de l'unit\'e est $\mathbb{Q}(\sqrt{-3})$. Par suite, le groupe $\mathcal{O}_{K}^{\times}/(\mathcal{O}_{K}^{\times})^p$ est nul, et la th\'eorie de Kummer sur $S$ nous dit alors que
$$
\begin{CD}
H^1(S,\mu_p) @>\sim>> \pic(\mathcal{O}_K)[p]. \\
\end{CD}
$$
D'autre part, le diagramme suivant est commutatif, \`a lignes exactes
$$
\begin{CD}
0 @>>> H^1(S,\mu_p) @>\sim>> \pic(\mathcal{O}_K)[p] @>>> 0 \\
@VVV @VVV @VVV \\
\mathcal{O}_{K'}^{\times}/(\mathcal{O}_{K'}^{\times})^p @>>> H^1(S',\mu_{p,S'})
@>>> \pic(\mathcal{O}_{K'})[p] @>>> 0 \\
\end{CD}
$$
et la fl\`eche verticale de droite est nulle. On en d\'eduit que la fl\`eche verticale du milieu est nulle, ce qu'on voulait.
\end{proof}

\begin{rmq}
Le choix le plus naturel est de prendre $K'$ contenu dans le corps de classes de Hilbert de $K$. Mais on peut \'egalement consid\'erer d'autres corps $K'$, \'eventuellement ramifi\'es sur $K$.
\end{rmq}

Nous sommes en mesure de d\'eduire du th\'eor\`eme \ref{essentiel} le r\'esultat suivant

\begin{thm}
\label{bigthm}
Soient
\begin{enumerate}
\item $p\geq 3$ un nombre premier,
\item $K$ un corps quadratique imaginaire tel que $\pic(\mathcal{O}_K)[p]\neq 0$ ; si $p=3$ on exige de plus que $p$ soit non ramifi\'e dans $K$,
\item $K'$ une extension de $K$ dans laquelle les \'el\'ements de $\pic(\mathcal{O}_K)[p]$ capitulent,
\item $E_{K'}$ une $K'$-courbe elliptique, ayant bonne r\'eduction au-dessus de $p$, telle que $E_{K'}(K')$ soit un groupe fini contenant un \'el\'ement d'ordre $p$.
\end{enumerate}
Alors il existe une $K$-vari\'et\'e ab\'elienne $A_K$ de dimension $[K':K]$, et une immersion ferm\'ee $\mu_p\rightarrow \A$ satisfaisant les hypoth\`eses de la proposition \ref{classinvariant}, telles que les groupes $B(S)$ et $\B^{\Lambda}(S)$ soient de torsion. Par cons\'equent, les homomorphismes $\psi$ et $\psibemol$ correspondants sont non nuls sur les points de torsion.
\end{thm}

\begin{proof}
Tout d'abord, $p$ est peu ramifi\'e dans $S$ (si $p>3$ cela d\'ecoule du fait que $[K:\mathbb{Q}]=2$, et si $p=3$ c'est vrai par hypoth\`ese). De plus, la fl\`eche de changement de base $H^1(S,\mu_p)\rightarrow H^1(S',\mu_{p,S'})$ est nulle d'apr\`es le lemme \ref{naif}. Il s'ensuit que les hypoth\`eses du th\'eor\`eme \ref{essentiel} sont satisfaites. La proposition \ref{classinvariant} s'applique alors, et implique la surjectivit\'e de $\psi$ et $\psibemol$ sur le groupe $\pic(\mathcal{O}_K)[p]$. Enfin, le groupe $B(S)$ est de m\^eme rang que le groupe $E_{K'}(K')$, et ce dernier est de torsion par hypoth\`ese.
\end{proof}

\begin{rmq}
Supposons que $E_{K'}$ provienne d'une $K$-courbe elliptique $E_K$, et que $E_K(K)$ contienne un point d'ordre $p$ (les seuls cas possibles sont $p=2$, $3$, $5$, $7$, $11$, $13$ d'apr\`es Kamienny \cite{kam}). Alors on peut prendre pour $K'$ le corps de classes de Hilbert, qui est non ramifi\'e sur $K$. Par suite, $A_K$ a bonne r\'eduction en les m\^emes places que $E_K$, en vertu du lemme \ref{lem16}. Si en outre les premiers de mauvaise r\'eduction de $E_K$ sont des id\'eaux principaux de $K$, alors en les inversant on trouve un contre-exemple \`a l'annulation de l'homomorphisme de classes pour les sch\'emas ab\'eliens.
\end{rmq}

\begin{exm}
Soit $K=\mathbb{Q}(\sqrt{-23})$, alors $\pic(\mathcal{O}_K)\simeq\mathbb{Z}/3\mathbb{Z}$ et $3$ est non ramifi\'e dans $K$. Soit $E$ la courbe elliptique (14A1 dans les tables \cite{cr}) d\'efinie par l'\'equation
$$
y^2+xy+y=x^3+4x-6.
$$
Nous avons $E(\mathbb{Q})\simeq\mathbb{Z}/6\mathbb{Z}$. Soit $H$ le corps de classes de Hilbert de $K$, un calcul avec \texttt{Magma} montre que $\rk(E(H))=0$. Les hypoth\`eses du th\'eor\`eme \ref{bigthm} sont alors satisfaites.
\end{exm}

\begin{exm}
Soit $K=\mathbb{Q}(\sqrt{-47})$, alors $\pic(\mathcal{O}_K)\simeq\mathbb{Z}/5\mathbb{Z}$. Soit $E$ la courbe elliptique (11A1 dans les tables \cite{cr}) d\'efinie par l'\'equation
$$
y^2+y=x^3-x^2-10x-20.
$$
Nous avons $E(\mathbb{Q})\simeq\mathbb{Z}/5\mathbb{Z}$. Soit $H$ le corps de classes de Hilbert de $K$, un calcul avec \texttt{Magma} montre que $\rk(E(H))=0$. Les hypoth\`eses du th\'eor\`eme \ref{bigthm} sont alors satisfaites. En outre, $(11)$ est un id\'eal premier de $K$, donc si l'on pose $U=\Spec(\mathcal{O}_K[11^{-1}])$, alors $\pic(U)\simeq\pic(\mathcal{O}_K)$, et $\A_U$ est un $U$-sch\'ema ab\'elien. Ainsi on obtient un exemple de sch\'ema ab\'elien de dimension $5$, muni d'une $5$-isog\'enie pour laquelle l'homomorphisme $\psi$ n'est pas nul sur les points de torsion.
\end{exm}

\noindent\textit{Remerciements}.
Je remercie chaleureusement Christian W\"uthrich pour son aide dans le calcul des exemples ci-dessus, et pour les discussions fructueuses que nous avons eues. Merci \'egalement \`a Sylvia Guibert pour ses suggestions qui m'ont permis d'am\'eliorer la structure de cet article.




\vskip 2cm

Jean Gillibert
\smallskip

The University of Manchester

Alan Turing Building

Oxford Road

Manchester M13 9PL

United Kingdom

\bigskip

\texttt{jean.gillibert@manchester.ac.uk}

\end{document}